\documentclass[12pt,twoside,final,reqno]{amsart}
\usepackage{amssymb,amsmath,amscd,amsfonts,amsthm,mathrsfs}
\usepackage{verbatim,paralist}
\usepackage{enumerate}
\RequirePackage{ifthen}
\usepackage[latin1]{inputenc}
\usepackage[T1]{fontenc}
\usepackage[english]{babel}
\usepackage{babel}
 \usepackage{tikz}
\usepackage{latexsym}
\usepackage{amsmath}
\usepackage{amsfonts}
\usepackage{amssymb}
\usepackage{latexsym}
\usepackage{color}
\usepackage{fancyhdr}
\newcommand\blue[1]{\textcolor{blue}{#1}}

\usepackage[latin1]{inputenc}

\newtheorem{theorem}{Theorem}[section]
\newtheorem{Proposition}[theorem]{Proposition}
\newtheorem{Lemma}[theorem]{Lemma}
\newtheorem{Remark}[theorem]{Remark}
\newtheorem{Definition}[theorem]{Definition}

\title {Browder $S$-resolvent equation in quaternionic setting}
\author{Hatem Baloudi}
\address{Hatem Baloudi, Department of Mathematics, Faculty of Sciences of Gafsa, University of Gafsa, 2112 Zarroug,
Tunisia}
\author{Aref JERIBI}
\address{Aref Jeribi, Department of Mathematics, Faculty of Sciences of Sfax, University of Sfax, 3000 Sfax,
Tunisia}
\author{Habib Zmouli}
\address{Habib Zmouli, Departement of Mathematics Faculty of Sciences of Sfax, University of Sfax, Route de Soukra km 3.5, B.P. 1171, 3000 Sfax, Tunisia.}

\begin{document}

\maketitle
\begin{abstract}
This paper is devoted to the study of the $S$-eigenvalue of finite type of a bounded right quaternionic linear operator acting in a right quaternionic Hilbert space. The study is based on the different properties of the Riesz projection associated with the connected part of the $S$-spectrum. Furthermore, we introduce the left and right Browder $S$-resolvent operators. Inspired by the $S$-resolvent equation, we give the Browder's $S$-resolvent equation in quaternionic setting.
\end{abstract}

\tableofcontents

\section{Introduction}
In the theory of complex Banach spaces, the search for the eigenvalues of finite type of linear operator aroused the interest and attracted the attention of many researchers, see for instance \cite{ABJK,BJ,Barnes,Bruce,Gohberg,J1,J2} and references therein. Sometimes this type of eigenvalue is known as a Riesz point. In doing so, one can develop the spectral theory of operators. An attractive characterization of eigenvalues of finite type by using Riesz projection is discussed and determined in \cite{Gohberg}. In particular, they show that this part of the spectrum is only the set of isolated point of the spectrum such that the corresponding Riesz projections are finite dimensional. Thus an extension of the usual resolvent is studied in \cite{Lutgen}. We refer to \cite{charfi,J1} for the applications to the Frobenius Schur factorizationfor $2\times 2$ Matrices and to the transport operators.
\vskip 0.2 cm

In the quaternionic setting, there has long been an apparent problem in defining the concept of spectrum of a quaternionic operator. In fact, the quaternionic multiplication is not commutative. This makes it possible to observe three types of Banach spaces: right, left and bilateral, according to the operation of the multiplication on the vectors.  It was only in $2006$ that F. Colombo and I. Sabbadini succeed in giving a new attractive and useful concept for the study of quaternionic operations, namely the $S$-spectrum. We refer to \cite[Section 1.2.1]{FJDP}, see also \cite{FJGanter} for the precise history and motivation of this new concept. Some years later, D. Alpay, F. Colombo and D.P. Kimsey in \cite{DFDP} gave the spectral theorem for the bounded and unbounded quaternionic operator related to the concept of $S$-spectrum. In the book \cite{NCFCBO}, the authors have studied and discussed the spectral theory for the Clifford operators. We refer to \cite{PF} for some results on operators perturbation, to \cite{FD1} for a version of functional calculus for bounded and unbounded normal operators on a Clifford module, and to \cite{slice} for the study and discussion of slice monogenic function of a Clifford variable.
\vskip 0.2 cm

The first aim of this article is to study the $S$-eigenvalue of finite type of a bounded right quaternionic linear operator acting in a right quaternionic Hilbert space. In fact, if $T\in \mathcal{B}(V^{R}_{\mathbb{H}})$ (the set of all right bounded operator) and $q\in\sigma_{S}(T)\setminus\mathbb{R}$ (where $\sigma_{S}(T)$ denote the $S$-spectrum of $T$), then $[q]:=\{hqh^{-1}:\ q\in\mathbb{H}^*\}\subset\sigma_{S}(T)$ since the $S$-spectrum of $T$ is axially symmetric. In particular, $q$ is never isolated in $\sigma_{S}(T)$. However, we can speak of an isolated $2$-sphere in $\sigma_{S}(T)$. In this case, we can associate to $T$ and $[q]$ with a Riesz projection. Under there circumstances, the $S$-eigenvalue of finite type will be considered as an isolated $2$-sphere with an associated Riesz projection $P_{[q]}$ of finite rank. Especially, if $T\in \mathcal{B}(V_{\mathbb{C}})$ (i.e. $T$ is a linear operator acting on complex Banach space) and $\lambda$ is in the complex spectrum of $T$, then $[\lambda]=\{\lambda\}$ and this gives the usual version of the Riesz point in the complex case. We turn to the understanding of the $S$-eigenvalue of finite type. To begin with, let $T\in \mathcal{B}(V^{R}_{\mathbb{H}})$ and $q\in\sigma_{d}^{S}(T)$ (the set of $S$-eigenvalues of finite type), we refer to Section \ref{S.1} for precise definition. The first result of this paper characterizes  the range of the Riesz projection $P_{[q]}$ associated with the $2$-sphere $[q]$ and the operator $T$. Next, thank to the $S$-spectral mapping theorem \cite[Theorem 4.2.1]{FJDP}, we show that if we perturb the pseudo-resolvent $\mathcal{Q}_{q}(T):=T^2-2{\rm Re}(q)T+\vert q\vert^2\mathbb{I}_{V^{R}_{\mathbb{H}}}$ by the Riesz projection $P_{[q]}$ we obtain an invertible operator. We end the first part of the paper with a discussion on the localization of the $S$-eigenvalue of finite type of sequence of quaternionic operators.
\vskip 0.2 cm

The second aim of this article is to determine the quaternionic version of Browder's resolvent equation. Let $T$ be a linear operator acting on a complex Banach space $V_{\mathbb{C}}$. The spectrum of $T$ will be denoted by $\sigma(T)$ and the Riesz point will be denoted by $\sigma_{d}(T)$ (the set of isolated point $\lambda\in\mathbb{C}$ in the spectrum such that the corresponding Riesz projection $P_{\{\lambda\}}$ are finite dimensional). For $\lambda,\ \mu\in(\mathbb{C}\setminus\sigma(T))\cup\sigma_{d}(T)$, the Browder's resolvent equation is given by
\begin{equation}\label{e.b}
R_{B}^{-1}(\lambda,T)-R_{B}^{-1}(\mu,T)=(\lambda-\mu)R_{B}(\lambda,T)R_{B}(\mu,T)+M_{T}(\lambda,\mu),
\end{equation}
where
\begin{align*}R_{B}^{-1}(\lambda,T)=(T-\lambda\mid_{P^{-1}_{\{\lambda\}}(\{0\})})^{-1}(I-P_{\{\lambda\}})+P_{\{\lambda\}}\end{align*}
and
\begin{align*}M_{T}(\lambda,\mu)=R_{B}^{-1}(\lambda,T)([T-(\lambda+1)]P_{\{\lambda\}}-[T-(\mu+1)]P_{\{\mu\}})R_{B}^{-1}(\mu,T).\end{align*}
\noindent We refer to \cite{Lutgen} for a brief discussion and for a full proof. We turn to the case quaternionic. Set $T\in \mathcal{B}(V^{R}_{\mathbb{H}})$ and $q\in\sigma_{d}^{S}(T)$. Let $P_{[q]}$ denote the corresponding Riesz projector with rang and kernel denoted by $R(P_{[q]})$ and $N(P_{[q]})$, respectively. Thanks to the Riesz decomposition theorem \cite[Theorem 6]{PSK} in quaternionic setting, we have
\begin{align*}\sigma_{S}(T\mid_{R(P_{[q]})})=[q]\mbox{ and }\sigma_{S}(T\mid_{N(P_{[q]})})=\sigma_{S}(T)\setminus[q].\end{align*}
\noindent In this way, we can define the left Browder $S-$resolvent operator
\begin{align*}S_{L,B}^{-1}(q,T):=-[\mathcal{Q}_{q}(T)\mid_{N(P_{[q]})}]^{-1}(T-\overline{q}\mathbb{I}_{V^{R}_{\mathbb{H}}})(\mathbb{I}_{V^{R}_{\mathbb{H}}}-P_{[q]})-P_{[q]}\end{align*}
and the right Browder $S-$resolvent operator
\begin{align*}
S_{R,B}^{-1}(q,T):=-(T-\overline{q}\mathbb{I}_{V^{R}_{\mathbb{H}}})[\mathcal{Q}_{q}(T)\mid_{N(P_{[q]})}]^{-1}(\mathbb{I}_{V^{R}_{\mathbb{H}}}-P_{[q]})-P_{[q]}.\end{align*}
Motivated by this, we obtain a generalization of the classical Browder's resolvent equation (\ref{e.b}). Precisely, we have

$S^{-1}_{R,B}(s,T)S^{-1}_{L,B}(p,T)\mathcal{Q}_{s}(p)$\begin{align*}
&=[S^{-1}_{R,B}(s,T)-S^{-1}_{L,B}(p,T)]p+\overline{s}[S^{-1}_{L,B}(p,T)-S^{-1}_{R,B}(s,T)]\\
&+[S^{-1}_{R,B}(s,T)(T-(p+1)\mathbb{I}_{V^{R}_{\mathbb{H}}})P_{[p]}-(T-(s+1)\mathbb{I}_{V^{R}_{\mathbb{H}}})P_{[s]} S^{-1}_{L,B}(p,T)]p\\
&+\overline{s}[(T-(s+1)\mathbb{I}_{V^{R}_{\mathbb{H}}})P_{[s]}S^{-1}_{L,B}(p,T)-S^{-1}_{R,B}(s,T)(T-(p+1)\mathbb{I}_{V^{R}_{\mathbb{H}}})P_{[p]}],
\end{align*}
\noindent where $p, s\in(\mathbb{H}\setminus\sigma_{S}(T))\cup\sigma_{d}^S(T)$ and $\mathcal{Q}_{s}(p)=p^2-2{\rm Re}(s)p+\vert s\vert^2$. The technique of the proof is inspired from the proof of \cite[Theorem 3.1.15]{FJDP}. It is remarkable that the Browder's resolvent equation extend \cite[Theorem 3.1.15]{FJDP} to $(\mathbb{H}\setminus\sigma_{S}(T))\cup\sigma_{d}^S(T)$. Indeed, if $q\in \mathbb{H}\setminus\sigma_{S}(T)$ with the convention $P_{[q]}=0$, then
\begin{align*}S_{L,B}^{-1}(q,T)=S_{L}^{-1}(q,T)=-(T^2-2{\rm Re}(q)T+\vert q\vert^2\mathbb{I}_{V^{R}_{\mathbb{H}}})^{-1}(T-\overline{q}\mathbb{I}_{V^{R}_{\mathbb{H}}})\end{align*}
\noindent and
\begin{align*}S_{R,B}^{-1}(q,T)=S_{R}^{-1}(q,T)=-(T-\overline{q}\mathbb{I}_{V^{R}_{\mathbb{H}}})(T^2-2{\rm Re}(q)T+\vert q\vert^2\mathbb{I}_{V^{R}_{\mathbb{H}}})^{-1}.\end{align*}
\vskip 0.1 cm

 As for the rest of this paper, it is structured as follows. The next Section is devoted to some basic notions of operator theory and slice functional calculus. In Section 3, we discuss some properties of the $S$-eigenvalue of finite type. Finally, in Section 4, we give and provide the Browder's $S$-resolvent equation in quaternionic setting.

\section{Mathematical preliminaries} In order to make the paper detailed, we collect some definitions and recall some results needed in the rest of the paper. We refer to \cite{SL,DFDP,FJDP,FJGanter} for surveys on the matter.

\subsection{Quaternions}

We denote by $\mathbb{H}$ the Hamiltonian skew field of quaternions with the standard basis $\{1, i, j, k\}$. Formally, we have
\begin{align*}\mathbb{H}=\Big\{q=x_0+x_1i+x_2j+x_3k:\ x_i\in\mathbb{R},\ i=0,1,2,3\Big\}.\end{align*}
The three imaginary units $i,j,k$ satisfy the relations
\begin{align*}i^2=j^2=k^2=ijk=-1,\ ij=-ji=k,\ ki=-ik=j,\ jk=-kj=i.\end{align*}
Let $q=x_0+x_1i+x_2j+x_3k\in\mathbb{H}$. The real part of $q$ is given by ${\rm Re}(q)=x_0$ and its imaginary part is defined as ${\rm Im}(q)=x_1i+x_2j+x_3k$, then the conjugate and the usual norm of $q$ are defined, respectively, by
\begin{align*}\overline{q}={\rm Re}(q)-{\rm Im}(q)\mbox{  and  }\vert q\vert=\sqrt{x_0^2+x_1^2+x_2^2+x_3^2}.\end{align*}
The unit sphere of purely imaginary quaternions is given by
\begin{align*}\mathbb{S}=\Big\{q\in\mathbb{H}:\ {\rm Re}(q)=0\mbox{ and }\overline{q}q=1\Big\}.\end{align*}
It is remarkable that $\mathbb{S}$ is a two-dimensional sphere in $\mathbb{R}^4$. If $q\in\mathbb{H}\setminus \mathbb{R}$, then
\begin{align*}q={\rm Re}(q)+I_q\vert {\rm Im}(q)\vert, \end{align*}
where $I_q=\frac{{\rm Im}(q)}{\vert{\rm Im}(q)\vert}\in\mathbb{S}$. In this way we can associated to $q$ a two-dimensional sphere defined by
\begin{align*}[q]={\rm Re}(q)+\mathbb{S}\vert {\rm Im}(q)\vert.\end{align*}
Note that $[q]$ has center at the real point ${\rm Re}(q)$ and has  radius $\vert{\rm Im}(q)\vert$. This sphere $[q]$ coincides with the set $\{hqh^{-1}:\ h\in \mathbb{H}^*\}$. We refer the reader to \cite{DFIS} for the full proof. For $I\in\mathbb{S}$, we set
\begin{align*}\mathbb{C}_I=\mathbb{R}+I\mathbb{R}.\end{align*}
In this case, we have
\begin{align*}\mathbb{H}=\bigcup_{I\in\mathbb{S}}\mathbb{C}_{I}.\end{align*}

\subsection{Operators acting on right quaternionic Hilbert space}
 Let $V^{R}_{\mathbb{H}}$ be a  right quaternionic Hilbert space and $\mathcal{O}=\{\phi_k:\ k\in\mathbb{N}\}$ be an orthonormal subset of $V^{R}_{\mathbb{H}}$. $\mathcal{O}$ is said to be Hilbert basis of $V^{R}_{\mathbb{H}}$ if  for every $\phi,\psi\in V_{\mathbb{H}}^{R}$, the series $\sum_{k\in\mathbb{N}}\langle \phi,\phi_{k}\rangle\langle \phi_{k},\psi\rangle$ converges absolutely and
\begin{align*} \langle \phi,\psi\rangle=\sum_{k\in\mathbb{N}}\langle \phi,\phi_{k}\rangle\langle \phi_{k},\psi\rangle.\end{align*}

\noindent A more detailed discussion about right Hilbert space can be found  in \cite{SL,FJDP}. We start by recalling the following classical result whose proof is analogous to that in the complex version, see \cite{KV}.

\begin{Proposition}  \label{artprop2.1}
Let $V_{\mathbb{H}}^{R}$ be a separable right quaternionic Hilbert space. We have the following three assertions
\begin{enumerate}
\item $V_{\mathbb{H}}^{R}$ admits at least one Hilbert basis.
\item If $X$ and $Y$ are two Hilbert basis of $V_{\mathbb{H}}^{R}$, then $\sharp X=\sharp Y,$ where $\sharp X$ denote the cardinal of $X$.
\item If $\mathcal{O}=\{\phi_k:\ k\in\mathbb{N}\}$ is  a Hilbert basis of $V_{\mathbb{H}}^{R}$, then every $\phi\in V_{\mathbb{H}}^{R}$ can be uniquely decomposed as follows
\begin{align*}\phi=\sum_{k\in\mathbb{N}}\phi_{k}\langle \phi_{k},\phi\rangle,\end{align*}
\noindent where the series $\sum_{k\in\mathbb{N}}\phi_{k}\langle \phi_{k},\phi\rangle$ converges absolutely in $V_{\mathbb{H}}^{R}$.
\end{enumerate}
\end{Proposition}

\noindent Let $\mathcal{O}=\{\phi_{k}:\ , k\in \mathbb{N}\}$ be a Hilbert basis of $V^{R}_{\mathbb{H}}.$ The left scalar multiplication on $V^{R}_{\mathbb{H}}$ induced by $\mathcal{O}$ is defined as the map
\begin{align*} \mathbb{H}\times V^{R}_{\mathbb{H}} & \longrightarrow V^{R}_{\mathbb{H}}\\
                                             (q, \phi)& \longrightarrow q\phi= \sum _{k \in \mathbb{N}} \phi_{k}q\langle\phi_{k}, \phi\rangle. \end{align*}

 A function $T:\ V^{R}_{\mathbb{H}}\longrightarrow V^{R}_{\mathbb{H}}$ is said to be quaternionic right linear if
\begin{align*}T(\phi+\psi q)=T(\phi)+T(\psi)q,\end{align*}
for all $\phi$, $\psi$ $\in V^{R}_{\mathbb{H}}$ and $q\in \mathbb{H}.$
We call a quaternionic right operator $T$ is bounded if
\begin{align*}\|T\|\displaystyle=\sup_{\phi\in V^{R}_{\mathbb{H}}\backslash\{0\}} \frac{\|T\phi\|}{\|\phi\|}< +\infty.\end{align*}
The set of all bounded right operators on $ V^{R}_{\mathbb{H}}$ is denoted by $\mathcal{B}(V^{R}_{\mathbb{H}})$ and the identity operator on $V^{R}_{\mathbb{H}}$ will be denoted by $\mathbb{I}_{V^{R}_{\mathbb{H}}}.$
If $T$ $\in \mathcal{B}(V^{R}_{\mathbb{H}})$, then  we write $N(T)$ and $R(T)$, respectively, for the null space and range of $T$. We set $$\alpha(T)=\dim N(T)\mbox{ and } \beta(T)={\rm codim \ } R(T).$$

\begin{Definition}\cite{BK,BK2} \label{artdef2.2}
Let $T\in \mathcal{B}(V^{R}_{\mathbb{H}})$, then
\begin{enumerate}
\item $T$ is a Fredholm operator if both $\alpha(T)$ and $\beta(T)$ are finite.
\item If $T$ is a Fredholm operator, then the index of $T$ is the number
\begin{align*}i(T)=\alpha(T)-\beta(T).\end{align*}
\item $T$ is a Weyl operator if $T$ is a Fredholm operator and $i(T)=0$.
\end{enumerate}
\end{Definition}
\noindent Let $\Phi(V^{R}_{\mathbb{H}})$ be denote the set of Fredholm operators and $\mathcal{W}(V^{R}_{\mathbb{H}})$ be denote the set of Weyl operators.

\begin{Definition}\cite{BK,BK2} Let $T\in \mathcal{B}(V^{R}_{\mathbb{H}})$.
\begin{enumerate}
\item $T$ is said a finite rank if $\dim R(T)<\infty$.
\item $T$ is said compact if $T$ maps bounded set into precompact sets.
\end{enumerate}
\end{Definition}

\noindent We denote by $\mathcal{K}(V^{R}_{\mathbb{H}})$ the set of all compact operators on $V^{R}_{\mathbb{H}}$. In the sequel of the paper, we equip $V^{R}_{\mathbb{H}}$ with a Hilbert basis $\mathcal{O}$. In this way, $\mathcal{B}(V^{R}_{\mathbb{H}})$ is a two-sided ideal quaternionic Banach algebras with respect to the two multiplications:
 \begin{align*}(qT)\phi=\displaystyle \sum_{\psi\in \mathcal{O}}\psi q\langle \psi,T\phi\rangle\mbox{ and }(Tq)\phi=\displaystyle\sum_{\psi\in \mathcal{O}}T(\psi)q\langle \psi,\phi\rangle.\end{align*} for all $\phi \in V^{R}_{\mathbb{H}}$.
In the next proposition we will recall some well-known properties of the compact and Fredholm-set, see \cite{BK,BK2}.

\begin{Proposition}~~\label{artprop2.4}
\begin{enumerate}
\item $\mathcal{K}(V^{R}_{\mathbb{H}})$ is a closed two-sided ideal of  $\mathcal{B}(V^{R}_{\mathbb{H}}).$
\item If $A\in\Phi(V^{R}_{\mathbb{H}})$ and $K\in \mathcal{K}(V^{R}_{\mathbb{H}})$, then $A+K\in \Phi(V^{R}_{\mathbb{H}})$ and $i(A+K)=i(A)$.
\end{enumerate}
\end{Proposition}

\subsection{The quaternionic functional calculus}
In this subsection, we recall some definitions and basis properties for the  Sabadini spectrum ($S$-spectrum), slice regular functions and Riesz projectors necessary for development of this manuscript. For more details see \cite{A:C,DFDP,NCFCBO,FJDP,FJGanter}.
\vskip 0.1 cm

For $T\in \mathcal{B}(V^{R}_{\mathbb{H}})$ and $q\in\mathbb{H}$, we define the associated operator $\mathcal{Q}_q(T):\ V^{R}_{\mathbb{H}}\longrightarrow V^{R}_{\mathbb{H}}$ by setting
\begin{align*}\mathcal{Q}_{q}(T):=T^2-2{\rm Re}(q)T+\vert q\vert^2\mathbb{I}_{V^{R}_{\mathbb{H}}}.\end{align*}

\begin{Definition}\label{def2.5}
Let $T\in \mathcal{B}(V^{R}_{\mathbb{H}})$.
\begin{enumerate}
\item The $S$-spectrum of $T$ is defined as
\begin{align*}\sigma_{S}(T)=\Big\{q\in\mathbb{H}:\ \mathcal{Q}_{q}(T)\mbox{ is not invertible in } \mathcal{B}(V^{R}_{\mathbb{H}})\Big\}.\end{align*}
\item We define the $S$-resolvent set of $T$ as
\begin{align*}\rho_{S}(T)=\mathbb{H}\setminus\sigma_{S}(T).\end{align*}
\item The point $S$-spectrum of $T$ is given by
\begin{align*}\sigma_{pS}(T)=\Big\{q\in\mathbb{H}:\ N(\mathcal{Q}_q(T))\neq\{0\}\Big\}.\end{align*}
\end{enumerate}
\end{Definition}

 \noindent The concept of $S$-spectrum is motivated by both the left Cauchy kernel series
 \begin{align*}\displaystyle\sum_{n=0}^{+\infty}T^nq^{-n-1}=-(T^2-2{\rm Re}(q)T+\vert q\vert^2\mathbb{I}_{V^{R}_{\mathbb{H}}})^{-1}(T-\overline{q}\mathbb{I}_{V^{R}_{\mathbb{H}}}),\ \vert q\vert >\parallel T\parallel\end{align*}
  and the right Cauchy kernel series
 \begin{align*}\displaystyle\sum_{n=0}^{+\infty}q^{-n-1}T^n=-(T-\overline{q}\mathbb{I}_{V^{R}_{\mathbb{H}}})(T^2-2{\rm Re}(q)T+\vert q\vert^2\mathbb{I}_{V^{R}_{\mathbb{H}}})^{-1},\ \vert q\vert >\parallel T\parallel.\end{align*}
 We refer to \cite{FJDP} for a full explanation.  Note that $\sigma_{S}(T)$ is a non-empty compact set, see \cite{FJGanter}. If $Tu=uq$ for some $u\in V^{R}_{\mathbb{H}}\backslash\{0\}$ and $q\in\mathbb{H}$, then $u$ is called eigenvector of $T$ with right eigenvalue $q$. We recall that the set of right eigenvalue coincides with point $S$-spectrum $\sigma_{S}(T)$, see \cite [Theorem 2.5]{CFSI}.

 \begin{Definition} \label {def2.6}
{\rm A set $\Omega\subset\mathbb{H}$ is called\\
\noindent $(i)$ axially symmetric if $\{hqh^{-1}:\ h\in\mathbb{H}\}\subset \Omega$ for any $q\in \Omega$ and\\
\noindent $(ii)$ a slice domain (or $s$-domain for short) if $\Omega$ is open, $\Omega\cap\mathbb{R}\neq\emptyset$ and $\Omega\cap\mathbb{C}_{I}$ is a domain in $\mathbb{C}_{I}$, for any $I\in\mathbb{S}$.}
\end{Definition}

\noindent Note that The $S$-spectrum $\sigma_{S}(T)$ and the $S$-resolvent $\rho_{S}(T)$ are axially symmetric, see \cite{FJDP}.

\begin{Definition}\cite[Definition 2.1.2]{FJDP} (Slice hyperholomorphic functions) Let $\Omega\subset\mathbb{H}$ be an axially symmetric open set and $f:\ \Omega\longrightarrow \mathbb{H}$ be a function. Set
\begin{align*}\Omega_{\mathbb{R}^2}:=\Big\{(u,v)\in\mathbb{R}^2:\ u+Iv\in\Omega\mbox{, for all }I\in\mathbb{S}\Big\}.\end{align*}
We say that $f$ is a left slice hyperholomorphic function if it is of the form
\begin{align*}f(q)=P(u,v)+I_{q}Q(u,v)\mbox{, for }q=u+I_{q}v\in\Omega\end{align*}
with $P$,$Q$ take value in $\mathbb{H}$ such that
\begin{equation}\label{r1}P(u,-v)=P(u,v),\ Q(u,-v)=-Q(u,v)\end{equation}
and satisfy the Cauchy Riemann equation
\begin{equation}\label{r2}
\displaystyle\frac{\partial P}{\partial u}-\frac{\partial Q}{\partial v}=0,\ \frac{\partial P}{\partial v}+\frac{\partial Q}{\partial u}=0.
\end{equation}
We say that $f$ is a right hyperholomorphic function if it is of the form
\begin{align*}f(q)=P(u,v)+Q(u,v)I_{q}\mbox{ for }q=u+I_{q}v\in\Omega\end{align*}
with $P$,$Q$ $:\ \Omega_{\mathbb{R}^2}\longrightarrow \mathbb{H}$ satisfy $(\ref{r1})$ and $(\ref{r2})$.\\
\noindent If $f$ is left or right with $P(u,v),\ Q(u,v)\in\mathbb{R}$ for all $(u,v)\in \Omega_{\mathbb{R}^2}$, then $f$ is said intrinsic function.
\end{Definition}
\vskip 0.2 cm
\noindent Let $\mathcal{SH}_{L}(\Omega)$ $($resp. $\mathcal{SH}_{R}(\Omega))$ be denote the set of left $($resp. right) slice hyperholomorphic functions on $\Omega$ and $\mathcal{N}(\Omega)$ be denotes the set of intrinsic functions. This class of functions  is a generalization of the set of holomorphic functions in the complex setting.

\begin{Definition} \label{def2.8}
{\rm Let $T\in \mathcal{B}(V_{\mathbb{H}}^{R})$ and $q\in\rho_{S}(T)$. The \emph{left $S$-resolvent operator} is defined by
\begin{align*}S_{L}^{-1}(q,T):=-(T^2-2{\rm Re}(q)T+\vert q\vert^2\mathbb{I}_{V^{R}_{\mathbb{H}}})^{-1}(T-\overline{q}\mathbb{I}_{V_{\mathbb{H}}^{R}}),\end{align*}
and the \emph{right $S$-resolvent operator} is given by
\begin{align*}S_{R}^{-1}(q,T):=-(T-\overline{q}\mathbb{I}_{V_{\mathbb{H}}^{R}})(T^2-2{\rm Re}(q)T+\vert q\vert^2\mathbb{I}_{V^{R}_{\mathbb{H}}})^{-1}.\end{align*}}
\end{Definition}

\begin{Proposition}\cite[Lemma 3.1.11]{FJDP} \label{prop2.9}
The left $S$-resolvent operator $q\longmapsto S_L^{-1}(q,T)$ is right slice hyperholomorphic and the right $S$-resolvent operator $q\longrightarrow S_{R}^{-1}(q,T)$ is left slice hyperholomorphic.
\end{Proposition}

Let $\mathcal{SH}_{L}(\sigma_{S}(T))$,\ $\mathcal{SH}_{R}(\sigma_{S}(T)))$ and  $\mathcal{N}(\sigma_{S}(T))$ be denote, respectively, the set of all  $f$ left, right and intrinsic slice hyperholomorphic functions $f$ such that $\sigma_{S}(T)\subset \mathcal{D}(f)$, where $\mathcal{D}(f)$ denote the domain of $f$.

\begin{Remark}\cite[Remark 3.2.4]{FJDP} \label{rem 2.10}
Let $f\in \mathcal{SH}_{L}(\sigma_{S}(T))\cup \mathcal{SH}_{R}(\sigma_{S}(T)))\cup \mathcal{N}(\sigma_{S}(T))$, then there exists a bounded slice Cauchy domain $\Omega$ such that
\begin{align*}\sigma_{S}(T)\subset \Omega\mbox{  and  }\overline{\Omega}\subset\mathcal{D}(f).\end{align*}
\end{Remark}

\noindent Now, we can give the version of the quaternionic functional calculus.

\begin{Definition}\cite[Definition 3.2.5]{FJDP} \label{def2.11}
Let $T\in \mathcal{B}(V_{\mathbb{H}}^{R})$. We define
\begin{equation}\label{e:4}f(T)=\displaystyle\frac{1}{2\pi}\int_{\partial(\Omega\cap\mathbb{C}_{I})}S_{L}^{-1}(q,T)dq_{I}f(q),\ \forall f\in \mathcal{SH}_{L}(\sigma_{S}(T)) \end{equation}
and
\begin{equation}\label{e:5}f(T)=\displaystyle\frac{1}{2\pi}\int_{\partial(\Omega\cap\mathbb{C}_{I})}f(q)dq_{I}S_{R}^{-1}(q,T),\ \forall f\in  \mathcal{SH}_{R}(\sigma_{S}(T))\end{equation}
where $dq_{I}=-dqI$ and $\Omega$ is a slice Cauchy domain as in the remark \ref{rem 2.10}.
\end{Definition}

\begin{theorem}$($Riesz's projectors$)$\cite[Theorem 4.1.5]{FJDP}
Let $T\in \mathcal{B}(V^{R}_{\mathbb{H}})$ and assume that $\sigma_{S}(T)= \sigma_{1} \cup\sigma_{2}$ with
\begin{align*}  dist(\sigma_{1},\sigma_{2}) > 0.\end{align*}
 Let $\mathcal{O}$ be  an open axially symmetric set with $\sigma_{1}\subset \mathcal{O}$ and $\overline{\mathcal{O}} \cap \sigma_{2}=\emptyset.$ We
define $\chi_{\sigma_{1}}(s)=1$  for $s \in \mathcal{O}$ and $\chi_{\sigma_{2}}(s)=0$ for $s\notin \mathcal{O},$\
Then, $\chi_{\sigma_{1}} \in\mathcal{N}(\sigma_{S}(T)), $ \
and
 \begin{align*}
 P_{\sigma_{1}}:=\chi_{\sigma_{1}}(T) = \frac{1}{2\pi}\int_{\partial(\mathcal{O}\cap\mathbb{C}_{I})} S^{-1}_{L}(q,T)dq_{I}.
 \end{align*}
Further, $P_{\sigma_{1}}$ is a continuous projection operator that commute with $T$ and $P_{\sigma_{1}}V^{R}_{\mathbb{H}}$ is a right linear subspace of $V^{R}_{\mathbb{H}}$  that is invariant under $T$.
\end{theorem}

\begin{theorem} \cite[Lemma 4.1.1]{FJDP}\label{artlem1} \label {artth2.13}
Let T $\in\mathcal{B}(V^{R}_{\mathbb{H}}),$ then
\begin{enumerate}
\item If $f$, $g$ $\in\mathcal{SH}_{L}(\sigma_{S}(T))$ and $q\in\mathbb{H},$ then
\begin{align*} (f+g)(T)=f(T)+g(T)   \ and  \ (fq)(T)=f(T)q.   \end{align*}
\item If $f$, $g$ $\in\mathcal{SH}_R(\sigma_{S}(T))$ and $q\in\mathbb{H},$ then
\begin{align*} (f+g)(T)=f(T)+g(T)   \ and  \ (qf)(T)=qf(T).   \end{align*}
\end{enumerate}
\end{theorem}

\begin{theorem}(The spectral mapping theorem)  \cite[Theorem 4.2.1]{FJDP}\label{artth2.14}\\
Let $T\in \mathcal{B}(V^{R}_{\mathbb{H}})$ and $f \in \mathcal{N}(\sigma_{S}(T)),$ then
\begin{align*}
\sigma_{S}(f(T))=f(\sigma_{S}(T)).
\end{align*}
\end{theorem}

\begin{Remark}\label{rem2.15}
Let $P_{\sigma}$ be a Riesz projector associated to the spectral set $\sigma,$ then
\begin{align*} qP_{\sigma}= P_{\sigma}q, \mbox{ for all } q\in\mathbb{H}. \end{align*}
In particular $R(P_{\sigma})$ is a left linear subspace of $V^{R}_{\mathbb{H}}.$ Indeed,
\begin{align*}
(qP_{\sigma})(T)=qP_{\sigma}(T)= q\chi_{\sigma}(T)=(\chi_{\sigma}q)(T)= \chi_{\sigma}(T)q = P_{\sigma}q.
\end{align*}
\end{Remark}

\section{Eigenvalue of finite type}\label{S.1}

\noindent Let
\begin{align*}\mathbb{R}_+^2:=\{(x,y)\in\mathbb{R}^{2}:\ y\in\mathbb{R}_+\}.\end{align*}

 We consider the following equation
 \begin{align*}
&
\Psi: \mathbb{H} \longrightarrow \mathbb{R}_+^2
\\
&\quad \quad q \longmapsto\ \ ({\rm Re}(q),\vert {\rm Im(q)}\vert).
\end{align*}
\noindent We refer to \cite{GN0} for more properties of $\Psi$. In particular, the author prove that $\Psi$ is continuous, open, closed and
\begin{align*} [\Omega]=\Psi^{-1}(\Psi(\Omega))\mbox{ for all }\Omega\subset\mathbb{H},\end{align*}
see \cite[Corollary 3.16 and Lemma 3.18]{GN0}.\\
Let $T \in \mathcal{B}(V_{\mathbb{H}}^R)$. A subset $\sigma\subset\sigma_{S}(T)$ is called an isolated part of $\sigma_{S}(T)$ if both $\sigma$ and $\sigma_{S}(T)\backslash\sigma$ are closed subsets of $\sigma_{S}(T)$. We start with the following result:

\begin{Proposition} \label{artprop3.1}
Let $T \in \mathcal{B}(V_{\mathbb{H}}^R)$ and $q\in\sigma_{S}(T)$, then $[q]$ is an isolated part of $\sigma_{S}(T)$ if, and only if, there exist $\varepsilon>0$ such that
\begin{equation}\label{e}[B(q,\varepsilon)]\cap\sigma_{S}(T)=[q],\end{equation}
where $B(q,\varepsilon)$denote the open boule of center $q$ and radius $\varepsilon.$
\end{Proposition}

\proof If $[q]$ is an isolated $2$-sphere of $\sigma_{S}(T)$, then $[q]$ is an open set of $\sigma_{S}(T)$. Let $U_q$ be an open set of $\mathbb{H}$ such that
\begin{align*}[q]=\sigma_{S}(T)\cap U_q.\end{align*}
 Since $q\in U_q$, then there exists $\varepsilon>0$ such that $B(q,\varepsilon)\subset U_q$. This implies that
 \begin{align*}[B(q,\varepsilon)]\cap\sigma_{S}(T)=[q].\end{align*}
 Indeed, assume that there exists $p\in [B(q,\varepsilon)]\cap\sigma_{S}(T)\backslash[q]$. So, $p\in [q']$ for some $q'\in U_{q}\backslash [q]$. Since $\sigma_S(T)$ is axially symmetric, then $q'\in\sigma_{S}(T)$, contradiction.\\
Conversely, if $(\ref{e})$ is satisfied, then $[q]$ is open in $\sigma_{S}(T)$. Since $[q]$, is closed, we deduce that $[q]$ is an isolated part of $\sigma_S(T)$. \qed
\vskip 0.2 cm

We recall that in a complex setting, the eigenvalue of finite type is introduced and studied in \cite{Gohberg}. In particular, the authors gave a characterization of this type of spectrum by using the Riesz projection. A version in the quaternionic case is introduced in \cite{BBj} in the following definition.

\begin{Definition}\label {artdef3.2}
Let {\rm $T\in\mathcal{B}(V_{\mathbb{H}}^{R})$. A point $q\in\sigma_{S}(T)$ is called a $S$-eigenvalue of finite type if $V_{\mathbb{H}}^{R}$ is a direct sum of $T$-invariant subspaces $V_{1,\mathbb{H}}^{R}$
and $V_{2,\mathbb{H}}^{R}$ such that
\vskip 0.1 cm

$(H1)$ $\dim(V_{1,\mathbb{H}}^{R})<\infty$,
\vskip 0.1 cm

$(H2)$ $\sigma_{S}(T\vert_{V_{1,\mathbb{H}}^{R}})\cap \sigma_{S}(T\vert_{V_{2,\mathbb{H}}^{R}})=\emptyset$,
\vskip 0.1 cm

$(H3)$ $\sigma_{S}(T\vert_{V_{1,\mathbb{H}}^{R}})=[q]$.}

\end{Definition}
\noindent We start by recalling the following decomposition theorem:

\begin{theorem}\label{j}\cite[Theorem 4.4]{KRE}
Let $T \in \mathcal{B}(V_{\mathbb{H}}^R)$. Suppose that $P_1$ is a projector in $\mathcal{B}(V_{\mathbb{H}}^R)$ commuting with $T$ and set $P_2=\mathbb{I}_{V_{\mathbb{H}}^R}-P_1$. Let $V_j=P_j(V_{\mathbb{H}}^R)$, $j=1,2$, and define the operators $T_j=TP_j=P_jT$. Denote by $\widetilde{T}_j:=T_j\mid_{V_{j}}$, $j=1,2$, then
\begin{align*}\sigma_{S}(T)=\sigma(\widetilde{T}_1)\cup\sigma(\widetilde{T}_2).\end{align*}
\end{theorem}

\begin{Remark} \label {artrem3.4}
Let $T \in \mathcal{B}(V_{\mathbb{H}}^R)$ and assume that $[q]$ is an isolated part of $\sigma_{S}(T)$. By using Theorem \ref{j} and \cite[Theorem 3.7.8]{FJGanter}, we have
\begin{align*}\sigma_{S}(T)=[q]\cup\sigma_{S}(T(\mathbb{I}_{V_{\mathbb{H}}^R}-P_{[q]}))\mbox{ and }[q]\cap\sigma_{S}(T(\mathbb{I}_{V_{\mathbb{H}}^R}-P_{[q]}))=\emptyset,\end{align*}
where $P_{[q]}$ is the Riesz projection related to $[q]$ and $T$.
\end{Remark}

\noindent We recall:
\begin{theorem}\cite[Theorem 3.10]{BBj}\label{t:2}
Let $T\in\mathcal{B}(V_{\mathbb{H}}^{R})$ and $[q]$ be an isolated part of $\sigma_{S}(T)$, then $q$ is a right eigenvalue of finite type if and only if  $\dim R(P_{[q]})<\infty$.
\end{theorem}
We turn to the pseudo $S$-resolvent operator
\begin{align*}\mathcal{Q}_{q}(T)^{-1}:=(T^2-2{\rm Re}(q)T+\vert q\vert^2\mathbb{I}_{V_{\mathbb{H}}^R})^{-1},\ q\in\rho_{S}(T).\end{align*}
As in complex case, one generalize this concept by using the Riesz projection. Let $\sigma_d^S(T)$ be denote the set of all $S$-eigenvalues of $T\in\mathcal{B}(V_{\mathbb{H}}^{R})$ of finite type. By using \cite[Theorem 3.7.8]{FJGanter}, we have
\begin{align*}\sigma_{S}(T\mid_{R(P_{[q]})})=[q]\mbox{ and }\sigma_{S}(T\mid_{N(P_{[q]})})=\sigma_{S}(T)\backslash[q]\end{align*}
for all $q\in\sigma_d^S(T)$. Thus, $q\in\rho_{S}(T\mid_{N(P_{[q]})})$ for every $q\in\sigma_d^S(T)$. This allows us to extend the pseudo $S$-resolvent operator. More precisely, set
\begin{align*}\rho_{B,S}(T):=\rho_{S}(T)\cup\sigma_d^S(T).\end{align*}
If $q\in \rho_{B,S}(T)$, then the operator
\begin{align*}P\mathcal{Q}_{q}(T):=(T^2-2{\rm Re}(q)T+\vert q\vert^{2} \mathbb{I}_{V_{\mathbb{H}}^R})(\mathbb{I}_{V_{\mathbb{H}}^R}-P_{[q]})+P_{[q]}\end{align*}
is invertible and its inverse is given by
\begin{align*}PR_{B,S}(q,T)=((T^2-2{\rm Re}(q)T+\vert q\vert^{2} \mathbb{I}_{V_{\mathbb{H}}^R})\mid_{N(P_{[q]})})^{-1}(\mathbb{I}_{V_{\mathbb{H}}^R}-P_{[q]})^{-1}+P_{[q]}.\end{align*}
Using this new concept, we prove the following result.

\begin{Proposition}\label{artprop3.6}
\noindent Let $T\in \mathcal{B}(V^{R}_{\mathbb{H}})$. If $q\in \sigma^{S}_{d}(T)$ and $x\in V^{R}_{\mathbb{H}}$, then
\begin{align*} x\in R(P_{[q]}) \mbox{ if and only if } \lim_{n\longrightarrow +\infty}\|\mathcal{Q}_{q}^{n}(T)x\|^{\frac{1}{n}}=0,\end{align*}
\noindent where $P_{[q]}$ is the Riesz projection related to $[q]$ and $T$.
\end{Proposition}
\proof
Using the Riesz decomposition \cite[Theorem 6]{PSK}, we have
\begin{align*}
\sigma_{S}(TP_{[q]}|_{R(P_{[q]})})=[q] \mbox{ and }\sigma_{S}(T(\mathbb{I}_{V_{\mathbb{H}}^R}-P_{[q]})|_{N(P_{[q]})})=\sigma_{S}(T)\backslash[q].
\end{align*}
\noindent Set
\begin{align*}\mathcal{Q}_{q}(X)=X^{2}-2Re(q)X+|q|^{2}.\end{align*}

 \noindent Observe that $X\longmapsto \mathcal{Q}_{q}(X)\in\mathcal{N}(\sigma_{S}(TP_{[q]}\mid_{R(P_{[q]})}))$. On the other hand, the polynomial $\mathcal{Q}_{q}$ vanishes exactly at $[q]$, see \cite[Lemma 4.2.3]{DFIS}. Now, by the $S$-spectral mapping theorem \ref {artth2.14}, we have
\begin{align*}
\sigma_{S}(\mathcal{Q}_{q}(T)P_{[q]}\mid_{R(P_{[q]})})=\mathcal{Q}_{q}(\sigma_{S}(TP_{[q]}\mid_{R(P_{[q]})}))=\mathcal{Q}_{q}([q])=\{0\}.
\end{align*}
\noindent In particular, for $0\neq x\in R(P_{[q]})$, we have
\begin{align*}
\lim_{n\longrightarrow +\infty}\|\mathcal{Q}_{q}^{n}(T)x\|^{\frac{1}{n}}
&=\lim_{n\longrightarrow +\infty}\|\mathcal{Q}_{q}^{n}(T)P_{[q]}\mid_{R(P_{[q]})}x\|^{\frac{1}{n}}\\
&\leq\lim_{n\longrightarrow +\infty}\|Q_{q}^{n}P_{[q]}\mid_{R(P_{[q]})}\|^{\frac{1}{n}}\|x\|^{\frac{1}{n}}\\
&=r_{S}(\mathcal{Q}_{q}(T)P_{[q]}\mid_{R(P_{[q]})})=0.
\end{align*}
\noindent Conversely, set $0\neq x\in V^{R}_{\mathbb{H}}$ such that $\lim_{n\longrightarrow +\infty}\|\mathcal{Q}_{q}^{n}(T)x\|^{\frac{1}{n}}=0$. Take
\begin{align*}
  x_{n}:=P\mathcal{Q}_{q}^n(T)x=(\mathcal{Q}_{q}(T))^{n}(\mathbb{I}_{V_{\mathbb{H}}^R}-P_{[q]})x+P_{[q]}x\blue{.}\end{align*}
It is clear that
\begin{align*}\|(\mathbb{I}_{V_{\mathbb{H}}^R}-P_{[q]})x_{n}\|^{\frac{1}{n}} \leq \|(\mathbb{I}_{V_{\mathbb{H}}^R}-P_{[q]})\| \|(\mathcal{Q}_{q}(T))^{n}x\|^{\frac{1}{n}}\blue{.}\end{align*}
In this way, we say that
\begin{align*}\lim_{n\longrightarrow+\infty} \|(I-P_{[q]})x_{n}\|^{\frac{1}{n}}=0.\end{align*}
Since $q\in \rho^{S}_{B}(T)$, then
\begin{align*}x=PR_{B,S}^n(q,T)x_n=(\mathcal{Q}_{q}(T)\mid_{N(P_{[q]})})^{-n}(\mathbb{I}_{V_{\mathbb{H}}^R}-P_{[q]})x_{n}+P_{[q]}x_{n}.\end{align*}
Finally, we obtain
\begin{align*}\|(\mathbb{I}_{V_{\mathbb{H}}^R}-P_{[q]})x\|^{\frac{1}{n}}\leq\|(\mathbb{I}_{V_{\mathbb{H}}^R}-P_{[q]})x_{n}\|^{\frac{1}{n}} \times \|(\mathcal{Q}_{q}(T)P_{[q]}|_{N(P_{[q]})})^{-1}\|.\end{align*}
This implies that $\lim_{n\longrightarrow+\infty}\|(I-P_{[q]})x\|^{\frac{1}{n}}=0$ and so $x\in R(P_{[q]})$.\qed
\\

\noindent We recall that $T$ is quasi-nilpotent if $\sigma_{S}(T)=\{0\}.$

\begin{Proposition}\label{artprop3.7}
Let $T\in \mathcal{B}(V^{R}_{\mathbb{H}})$ and $[q]$ be an isolated $2$-sphere of $\sigma_{S}(T).$ Then,
\begin{enumerate}
\item $\mathcal{Q}_{q}(T)+P_{[q]}$ and $\mathcal{Q}_{q}(T)+[2T+(1-2Re(q))\mathbb{I}_{V^{R}_{\mathbb{H}}}]P_{[q]}$ are invertible,
\item $\mathcal{Q}_{q}(T)P_{[q]}$ and $\mathcal{Q}_{q}(T)[2T+(1-2Re(q))\mathbb{I}_{V^{R}_{\mathbb{H}}}]P_{[q]}$ are quasi-nilpotent.
\end{enumerate}
\end{Proposition}
\proof
\noindent $(1)$ Since $[q]$ is an isolated $2$-sphere of $\sigma_{S}(T)$, then there exist $\varepsilon>0$ such that $[\overline{B}(q,\varepsilon)]\cap\sigma_{S}(T)=[q]$. Set \begin{align*}U:=[B(q,\varepsilon)]\mbox{ and }V:=\mathbb{H}\backslash [\overline{B}(q,\varepsilon)].\end{align*}
Observe that $U$ and $V$ are two axially symmetric open sets,
\begin{align*}U\cap V=\emptyset,\ [q]\subset U \mbox{ and }\sigma_{S}(T)\setminus [q]\subset V.\end{align*}

 Let us define the functions
\begin{align*}g(p):=\left\{ \begin{array}{rl}
\noindent 1& \mbox{for } p\in U,
\\~~
\\
 \displaystyle
    0&
  \mbox{for } p \in V.
\end{array} \right.\end{align*}
and
\begin{align*} h(p):=p^{2}-2Re(q)p+|q|^{2},\ p\in\mathbb{H}.
\end{align*}
\noindent Then
\begin{align*}g(T)=P_{[q]}\mbox{ and }h(T)=\mathcal{Q}_{q}(T).\end{align*}
\noindent Recall that $h(p)=0$ if, and only, if $p\in[q]$, see \cite[Lemma 4.2.3]{DFIS}. So, $(g+h)(p)\neq 0$ for all $p\in\sigma_{S}(T)$. Indeed,
if $p\in [q]\subseteq U$, then $g(p)=1$ and $ h(p)= 0$. If $ p\in \sigma_{S}(T)\backslash [q]$, then $g(p)=0$ and \  $h(p)\neq 0.$ Now, by using the algebraic properties of the quaternionic functional calculus \cite[Lemma 4.1.1]{FJDP}, we have
\begin{align*}\mathcal{Q}_{q}(T)+P_{[q]}=(g+h)(T).\end{align*}
Finally, since $g+h\in\mathcal{N}(\sigma_{S}(T)),$ then thanks to the $S$-spectral mapping theorem \ref{artth2.14}, we conclude that $\mathcal{Q}_{q}(T)+P_{[q]}$ is an invertible operator.\\
\noindent We turn to the operator $\mathcal{Q}_{q}(T)+[2T+(1-2Re(q))\mathbb{I}_{V^{R}_{\mathbb{H}}}]P_{[q]}$. We consider the function

\begin{align*}k(p):=\left\{ \begin{array}{rl}
\noindent 2p+1-2Re(q)\quad \quad &\mbox{if } p\in U,
\\~~
\\
 0\quad \quad\quad \quad \quad \quad \quad&\mbox{if } p\in V.
\end{array} \right.
\end{align*}
A similar argument as before, we have
\vskip 0.1 cm

\begin{align*}k(T)=[2T+(1-2Re(q))\mathbb{I}_{V^{R}_{\mathbb{H}}}]P_{[q]}\end{align*}
 and
 \begin{align*}(h+k)(T)=\mathcal{Q}_{q}(T)+[2T+(1-2Re(q))\mathbb{I}_{V^{R}_{\mathbb{H}}}]P_{[q]}.\end{align*}
 \vskip 0.1 cm
\noindent On the other hand, we have $(k+h)(p)\neq 0$ for all $p\in \sigma_{S}(T)$. Indeed, if $p={\rm Re}(q)+I_{p}\vert {\rm Im}(q)\vert$ for some $I_{p}\in\mathbb{S}$ (i.e, $p\in [q]$), then $h(p)=0$ (by using \cite[Lemma 4.2.3]{DFIS}) and
\begin{align*}k(p):
&=1-2Re(q)+2p\\
&=1+2I_{q}\vert {\rm Im}(q)\vert\neq 0.\end{align*}
 \noindent Now, if $p\in\sigma(T)\backslash [q]$, then we have easily $h(p)\neq \{0\}$ and $k(p)=0$. Finally, we can conclude that  $\mathcal{Q}_{q}(T)+[1-2Re(T)+2T]P_{[q]}$ is an invertible operator.\\
 \vskip 0.1 cm
 \noindent $(2)$ Since $(hg)(p)=h(p)g(p)=0$ and $(hk)(p)=0$ for all $p\in\sigma_{S}(T)$, then by using \cite[Lemma 3.2.8]{FJDP} and the $S$-spectral mapping theorem \ref{artth2.14}, we have
\begin{align*}\sigma_{S}(Q(T)P_{[q]})=\sigma_{S}(h(T)g(T))=\{0\}\end{align*}
and
\begin{align*}\sigma_{S}(\mathcal{Q}_{q}(T)[2T+(1-2Re(q))\mathbb{I}_{V^{R}_{\mathbb{H}}}]P_{[q]})=\sigma_{S}(h(T)k(T))=\{0\}.\end{align*}
This completes the proof.\qed

\begin{Proposition} \label{artprop3.8}
Let $T\in \mathcal{B}(V^{R}_{\mathbb{H}})$ and $[q]$ be an isolated $2$-sphere of $\sigma_{S}(T).$ Then,
\begin{align*}q\in\sigma_d^S(T) \mbox{ if and only if } \mathcal{Q}_{q}(T)\in \mathcal{W}(V^{R}_{\mathbb{H}}).\end{align*}
\end{Proposition}
\proof
Suppose that $q\in\sigma_d^S(T)$. By using the previous proposition \ref{artprop3.7}, we have $\mathcal{Q}_{q}(T)+P_{[q]}$ is an invertible operator. In this way, we see that $\mathcal{Q}_{q}(T)+P_{[q]}\in \mathcal{W}(V^{R}_{\mathbb{H}})$. Since $P_{[q]}\in\mathcal{K}(V^{R}_{\mathbb{H}})$, then thanks to proposition \ref{artprop2.4} we deduce that $\mathcal{Q}_{q}(T)\in \mathcal{W}(V^{R}_{\mathbb{H}})$.
\vskip 0.1 cm

Conversely, Let $\varepsilon>0$ such that $\sigma_{S}(T)\cap[B(q,\varepsilon)]=[q]$, then the Riesz projection $P_{[q]}$ associated with $T$ and $q$ is given by
\begin{align*}P_{[q]}
&=\displaystyle\frac{-1}{2\pi}\int_{\partial([B(q,\varepsilon)]\cap \mathbb{C}_{I})}\mathcal{Q}_{p}^{-1}(T)(T-\overline{p}\mathbb{I}_{V^{R}_{\mathbb{H}}})dp_{I}.\end{align*}
\noindent Now, denote by $\pi$ the natural quotient map into the Calkin algebra $\mathcal{C}(V_{\mathbb{H}^R})=\mathcal{B}(V_{\mathbb{H}^R})/\mathcal{K}(V_{\mathbb{H}^R})$, then $q\in\rho_{S}(\pi(T))$ and
\begin{align*}\pi(P_{[q]})=\displaystyle\frac{-1}{2\pi}\int_{\partial([B(q,\varepsilon)]\cap \mathbb{C}_{I})}\mathcal{Q}_{p}^{-1}(\pi(T))(\pi(T)-\overline{p}\mathbb{I}_{\mathcal{C}(V_{\mathbb{H}^R})})dp_{I}=0.\end{align*}
So, $P_{[q]}\in\mathcal{K}(V^{R}_{\mathbb{H}})$. This implies that $\mathbb{I}_{R(P_{[q]})}:\ R(P_{[q]})\longrightarrow R(P_{[q]})$ is compact, we deduce that $\dim(R(P_{[q]}))<\infty$. \qed

\begin{Remark}\label{artrem3.9}
As in the complex setting, we have if $T\in \mathcal{B}(V^{R}_{\mathbb{H}})$ is invertible and $N$ is a

 nilpotent operator that commute with $T$, then $T+N$ is also invertible. Indeed, let $m\in\mathbb{N}^{*}$ such that $N^m=0$. Then, $\mathbb{I}_{V^{R}_{\mathbb{H}}}+N$ is invertible and its inverse is given by
\begin{align*}(I+N)^{-1}=\displaystyle\sum_{k=0}^{m-1}(-1)^kN^k.\end{align*}
\noindent Since $TN=NT$, then $T^{-1}N$ is nilpotent. In this way, we see that $T+N=T(\mathbb{I}_{V^{R}_{\mathbb{H}}}+T^{-1}N)$ is invertible.
\end{Remark}

\begin{theorem}\label{s.}
Let $T_{n}$ and $T$ be belong to $\mathcal{B}(V_{\mathbb{H}}^{R})$ with $n\in \mathbb{N}$ and $\|T_{n}-T\| \longrightarrow 0.$ We suppose that $0\in \sigma_d^S(T).$
For an axially symmetric $V \subset \mathbb{H},$ we set
\begin{align*}E_{T}^{V \bigcap \sigma_{S}(T)} = (V \cap \sigma_{S}(T))/\cong ,\end{align*}
 \noindent where $p\cong q$ if, and only, if $p\in [q]$, then there exist $N \in \mathbb{N}$ and an open axially symmetric $ V_{0}\subset \mathbb{H}$  such that $\sharp E_{T_{n}}^{V_{0} \cap \sigma_{S}(T_{n})}<\infty$
and $V_{0}\bigcap \sigma_{S}(T_{n}) \subset \sigma_d^S(T_{n})$ for all $n\geq N.$
\end{theorem}

To prove this theorem, we first need to show the following results.

\begin{Lemma} \label{artlem1}
Let  $T$ and  $S$  $\in Inv(\mathcal{B}(V^{R}_{\mathbb{H}}))\ (i.e.\ 0\in\rho_S(T)\cap\rho_S(S)).$ We assume that $\|T-S\| \leq \frac{1}{2}\|S^{-1}\|^{-1},$ then
\begin{align*} \|T^{-1}-S^{-1}\|\leq 2 \|S^{-1}\|^{2} \|T-S\|.\end{align*}
\end{Lemma}
\proof
The proof is exactly similar to the proof of  \cite[Lemma 5, p.11]{BDJ} in the complex setting.

\begin{Definition}\cite[Definition 15, p.25]{BDJ} \label{artdef3.12}
Let $X$, $Y$ be  two topological spaces and let $\phi$ be a function defined on the space X and whose values are subsets of the space Y. The mapping $\phi$ is upper semi-continuous on $x_{0}$ if for each neighborhood $V_{\phi(x_{0})}$ of $\phi(x_{0})$, there exist a neighborhood $U_{x_{0}}$ of $x_{0}$ such that $$\phi(x)\subset V_{\phi(x_{0})}, \  x\in U_{x_{0}}.$$
$\phi$ is said to be upper semi-continuous if x is a point of upper semi-continuity for $\phi$ for each $x \in X.$
\end{Definition}

\begin{Lemma}\cite[Lemma 16, p.25]{BDJ}\label{artlem2}
Let $X$, $Y$ be metric spaces, let $Y$ be compact and let $\phi$ be a mapping of $X$ into the closed subsets of $Y$, then $\phi$ is upper semi-continuous if and only if the following conditions holds
 \begin{align*}x_{n} \in X,\ y_{n} \in \phi(x_{n}),\
x=\lim_{n\longrightarrow +\infty} x_{n},\  y=\lim_{n\longrightarrow+\infty}y_{n}\ \Longrightarrow y \in \phi(x).\end{align*}
\end{Lemma}

\begin{Proposition}\label{artprop1}
Let  $\psi _{\mathcal{B}(V^{R}_{\mathbb{H}})}$:  $T\longrightarrow \sigma_{S}(T)$ be the function defined on the space $\mathcal{B}(V^{R}_{\mathbb{H}})$ and whose values are in the compact subset of $\mathbb{H}.$ Then, $\psi _{\mathcal{B}(V^{R}_{\mathbb{H}})}$ is upper semi-continuous.
\end{Proposition}

\proof
Let $A \in \mathcal{B}(V^{R}_{\mathbb{H}})$ and $R>\|A\|.$ We set
\begin{align*}X_{A}:= \Big\{ B\in \mathcal{B}(V^{R}_{\mathbb{H}}): \|B\| \leq R \Big\}.\end{align*}
and
\begin{align*}Y_{A}:=\Big\{q \in \mathbb{H}: |q|\leq R\Big\}.\end{align*}
First, if
$\psi_{\mathcal{B}(V^{R}_{\mathbb{H}})}\mid_{X_A}: X_A \longrightarrow \{\mbox{ the compact subset of } Y_{A}\}$ is upper semi-continuous, then $\psi_{B(V^{R}_{\mathbb{H}})}$ is upper semi-continuous on $A$.
Secondly, let $(T_{n})_{n}$ be a sequence of operators in $X_A$, $q_n\in \sigma_{S}(T_{n})$, $$\lim_{n \longrightarrow +\infty}\|T_{n}-T\|=0 \mbox{  and  } \lim_{n \longrightarrow +\infty} |q_{n}-q|=0$$
We have to show that $q\in \sigma_{S}(T).$ Indeed, we assume that $\mathcal{Q}_{q}(T) \in Inv(\mathcal{B}(V^{R}_{\mathbb{H}})),$ then $$\mathcal{Q}_{q}(T)=\lim_{n\longrightarrow+\infty}\mathcal{Q}_{q_{n}}(T_{n}).$$
In fact,
 \begin{align*}\|\mathcal{Q}_{q}(T)-\mathcal{Q}_{q_n}(T_n)\| \leq \|T_{n}^{2}-T^{2}\|+\|2Re(q_{n})T_{n}-2Re(q)T\|+ ||q_{n}|^{2}-|q|^{2}|.\end{align*}
In this way, we see that
\begin{align*}\|\mathcal{Q}_{q}(T)-\mathcal{Q}_{q_{n}}(T_{n})\|\longrightarrow 0 \mbox{ since }  T \longrightarrow T^{2} \mbox { is continuous }.\end{align*}
Let $\varepsilon > 0$ be such that $B(\mathcal{Q}_{q}(T),\varepsilon) \subset Inv(\mathcal{B}(V^{R}_{\mathbb{H}})).$ Then, there exist $N_{q} \in \mathbb{N}$ such that $\mathcal{Q}_{q_{n}}(T_{n}) \in B(\mathcal{Q}_{q}(T),\varepsilon)$ for all $n \geq N_{q}.$ This is a contradiction since $q_n\in \sigma_{S}(T_n)$ for all $n$ $\in \mathbb{N}.$ \qed

\begin{Lemma} \label{artlem3}
Let $P$ and $Q$ be two projections in $\mathcal{B}(V^{R}_{\mathbb{H}}).$ We assume that $\|P-Q\| < 1$, then
\begin{enumerate}
\item $R(P)\cong R(Q).$
\item The operator $T=QP+(\mathbb{I}_{V^{R}_{\mathbb{H}}}-Q)(\mathbb{I}_{V^{R}_{\mathbb{H}}}-P)$ is bijective.
\item $T(R(P))\subset R(Q)$ and $T(N(P))\subset N(Q).$
\end{enumerate}
\end{Lemma}
\proof
The proof is exactly similar to the proof of \cite[Theorem 12.4]{DLA} in the complex setting.
\qed
\begin{Lemma}\label{h.}\cite[Lemma 3.1.3]{FJGanter}
Let $T\in \mathcal{B}(V^{R}_{\mathbb{H}}). $ The functions $q \longrightarrow \mathcal{Q}_{q}(T)^{-1}$ and $q \longrightarrow T\mathcal{Q}_{q}(T)^{-1}$ which are defined on $\rho_{S}(T)$ and take values in $\mathcal{B}(V^{R}_{\mathbb{H}})$ are continuous.
\end{Lemma}
\vskip 0.1 cm

\noindent {\bf Proof of Theorem \ref{s.}.}
Let $\varepsilon$ $\in ]0, 1[$ and $\mathcal{O}\subset \mathbb{H}$ be an open axially symmetric subset with $\mathcal{O} \cap \overline{B}(0,\varepsilon)=\emptyset$ and $\sigma_{S}(T)\backslash \{0\} \subset \mathcal{O}.$ By using Lemma \ref{h.} for all I $\in \mathbb{S},$ there is $M_{I} \geq 1$ such
\begin{align*}\sup_{q \in \partial (\mathbb{C}_{I} \bigcap B(0,\varepsilon))}\|\mathcal{Q}_{q}(T)^{-1}\| \leq M_{I}.\end{align*}
On the other hand, $\mathcal{O} \cup B(0,\varepsilon)$ is a neighborhood of $\sigma_{S}(T).$
By using Proposition \ref{artprop1}, there exist $N_\varepsilon > 0$ such that
\begin{align*}\sigma_{S}(T_{k})\subset B(0,\varepsilon)\cup \mathcal{O}\end{align*}

for all $k \geq N_{\varepsilon}$. We choose $N_{\varepsilon}$ large enough such that
\begin{align*}\|T^{2}_{k}-T^{2}+2Re(q)(T_{k}-T)\| \|T-\overline{q}\mathbb{I}_{V^{R}_{\mathbb{H}}}\|+\|T-T_{k}\| \leq \frac{1}{4M^{2}_{I}},\end{align*}
for all $q \in \partial (\mathbb{C}_{I} \bigcap B(0,\varepsilon))$.

In view of Lemma \ref{artlem1}, we have
\begin{align*}\|\mathcal{Q}_{q}(T)^{-1}-\mathcal{Q}_{q}(T_{k})^{-1}\| \leq 2M^{2}_{I} \|T^{2}_{k}-T^{2}+2Re(q)(T_{k}-T)\|,\end{align*}
$q \in \partial (\mathbb{C}_{I} \bigcap B(0,\varepsilon))$.
In this way, we see that
\begin{align*}
\|S_{L}^{-1}(q,T)-S_{L}^{-1}(q,T_{k})\| &\leq M_{I} \|T-T_{k}\|+\|\mathcal{Q}_{q}^{-1}(T)-\mathcal{Q}_{q}^{-1}(T_{k})\| \|T-\overline{q}\mathbb{I}_{V^{R}_{\mathbb{H}}}\|
\|\\
&\leq 2M_I^2\times\frac{1}{4M_{I}^2}=\frac{1}{2}<1.
\end{align*}
Let $P_{0}$ be the Riesz projection associated to $0$ and $T$. Set
$$P_{\sigma_{N_{\varepsilon}}^{k}}:=\frac{1}{2\pi}\int_{\partial(B(0,\varepsilon)\bigcap\mathbb{C_{I}})} S_{L}^{-1}(s,T_{k}) ds_{I},$$
where $\sigma_{N_{\varepsilon}}^{k}:=D(0,\varepsilon)\bigcap \sigma_{S}(T_{k}).$

By using Lemma \ref{artlem3}, we have
$$R(P_{0})\cong R(P_{\sigma^{k}_{N_{\varepsilon}}}).$$
In particular, $\dim R(P_{\sigma^{k}_{N_{\varepsilon}}}) <\infty$ for all $k\geq N_{\varepsilon}.$
Applying \cite[theorem 3.17]{BBj}, we have
\begin{align*}\sharp E_{T_{k}}^{\sigma^{k}_{N_{\varepsilon}}} <\infty,\end{align*}

 $q\in \sigma_d^S(T_{k})$ for all $k\geq N_{\varepsilon}$ and $q \in \sigma^{k}_{N_{\varepsilon}}.$ \qed

\section{Browder $S$-resolvent equation in quaternionic setting}
Let $[q]$  be an isolate  $2$-sphere  of $\sigma_{S}(T).$In view of spectral decomposition theorem, we have
\begin{align*} \sigma_{S}(T|_{R(P_{[q]})})=[q]  \ \mbox{ and }  \
 \sigma_{S}(T|_{N(P_{[q]})})= \sigma_{S}(T)\backslash [q]. \end{align*}
Let $T \in\mathcal{B}(V^{R}_{\mathbb{H}}).$ For $S\in\rho^{S}_{B}(T):=\rho_{S}(T)\cup\sigma^{S}_{d}(T).$\
We define the left Browder $S$-resolvent operator as
\begin{align*} S^{-1}_{L,B}(s,T)= -[Q_{s}(T)|_{N(P_{[s]})}]^{-1}(T-\overline{s}\mathbb{I_{V^{R}_{\mathbb{H}}}})(\mathbb{I_{V^{R}_{\mathbb{H}}}}-P_{[s]})-P_{[s]}.\end{align*}
and the right Browder $S$-resolvent operator as
\begin{align*}S^{-1}_{R,B}(s,T)=-(T-\overline{s}\mathbb{I_{V^{R}_{\mathbb{H}}}})[Q_{s}(T)|_{N(P_{[s]})}]^{-1}(\mathbb{I_{V^{R}_{\mathbb{H}}}}-P_{[s]})-P_{[s]}.\end{align*}

\begin{Remark}\label{rem4.1}
The Browder $S$-resolvent operator extend the $S$-resolvent operator to $\rho_{S}(T)\cup\sigma_d^S(T)$. Indeed, if $q\in\rho_{S}(T)$ with the convention $P_{[q]}=0$, we have
\begin{align*}S^{-1}_{R,B}(q,T)=S^{-1}_{R}(q,T)  \  \mbox{ and } \  S^{-1}_{L,B}(q,T)=S^{-1}_{L}(q,T).\end{align*}
\end{Remark}
\begin{theorem}\label{th4.2}
Let $T\in\mathcal{B}(V^{R}_{\mathbb{H}})$ and $q\in\rho^{S}_{B}(T):=\sigma^{S}_{d}(T)\cup\rho_{S}(T).$
Then, the left Browder $S$-resolvent operator satisfy the left Browder $S$-resolvent equation
\begin{align*} S^{-1}_{L,B}(q,T)(\mathbb{I_{V^{R}_{\mathbb{H}}}}-P_{[q]})q-T(\mathbb{I_{V^{R}_{\mathbb{H}}}}-P_{[q]})S^{-1}_{L,B}(q,T)+P_{[q]}=\mathbb{I_{V^{R}_{\mathbb{H}}}}. \end{align*}
and the right Browder $S$-resolvent operator satisfy the right Browder $S$-resolvent equation
\begin{align*} q(\mathbb{I_{V^{R}_{\mathbb{H}}}}-P_{[q]})S^{-1}_{R,B}(q,T)-S^{-1}_{R,B}(q,T)(\mathbb{I_{V^{R}_{\mathbb{H}}}}-P_{[q]})T+P_{[q]}=\mathbb{I_{V^{R}_{\mathbb{H}}}}. \end{align*}
\end{theorem}
\proof
Let $q\in\rho^{S}_{B}(T).$ It is clear that,
\begin{align*}(\mathbb{I_{V^{R}_{\mathbb{H}}}}-P_{[q]})q=q(\mathbb{I_{V^{R}_{\mathbb{H}}}}-P_{[q]}),\ P(\mathbb{I_{V^{R}_{\mathbb{H}}}}-P_{[q]})=(\mathbb{I_{V^{R}_{\mathbb{H}}}}-P_{[q]})P_{[q]}\end{align*}
 and
 \begin{align*} T (\mathcal{Q}_{q }(T)|_{N(P_{[q]})})^{-1}= (\mathcal{Q}_{q}(T)|_{N(P_{[q]})})^{-1} T\mid_{N(P_{[q]})}. \end{align*}
We obtain
\begin{align*}
&S^{-1}_{L,B}(q,T)(\mathbb{I_{V^{R}_{\mathbb{H}}}}-P_{[q]})q-T(\mathbb{I_{V^{R}_{\mathbb{H}}}}-P_{[q]})S^{-1}_{L,B}(q,T)\\
&=-[\mathcal{Q}_{q}(T)|_{N(P_{[q]})}]^{-1}(Tq-\vert q\vert^{2}\mathbb{I_{V^{R}_{\mathbb{H}}}})(\mathbb{I_{V^{R}_{\mathbb{H}}}}-P_{[q]}) +[\mathcal{Q}_{q}(T)|_{N(P_{[q]})}]^{-1}(T^{2}-T\overline{q})(\mathbb{I_{V^{R}_{\mathbb{H}}}}-P_{[q]})\\
&=[\mathcal{Q}_{q}(T)|_{N(P_{[q]})}]^{-1}\mathcal{Q}_{q}(T)|_{N(P_{[q]})} (\mathbb{I_{V^{R}_{\mathbb{H}}}}-P_{[q]}) \\
&=\mathbb{I_{V^{R}_{\mathbb{H}}}}-P_{[q]}.
\end{align*}
The right $S$-resolvent equation follows by similar computations.\qed

\begin{Remark} \label{rem4.3}
\begin{enumerate}
\item The left and the right $S$-resolvent equation implies,
\begin{align*}
S_{L,B}^{-1}(q,T)q-TS_{L,B}^{-1}(q,T)-(T-(q+1)\mathbb{I_{V^{R}_{\mathbb{H}}}})P_{[q]}=\mathbb{I_{V^{R}_{\mathbb{H}}}} \end{align*}
and
\begin{align*}qS_{R,B}^{-1}(q,T)-S_{R,B}^{-1}(q,T)T-(T-(q+1)\mathbb{I_{V^{R}_{\mathbb{H}}}})P_{[q]}=\mathbb{I_{V^{R}_{\mathbb{H}}}}.
\end{align*}
\item If $q\in\rho_{S}(T)$, then $P_{[q]}=0$. In this case, we obtain the two equations in \cite[Theorem 3.1.14]{FJDP}:
\begin{align*}
S_{L}^{-1}(q,T)q-TS_{L}^{-1}(q,T)=\mathbb{I_{V^{R}_{\mathbb{H}}}}
\end{align*}
and
\begin{align*}
qS_{R}^{-1}(q,T)-S_{R}^{-1}(q,T)T=\mathbb{I_{V^{R}_{\mathbb{H}}}}.
\end{align*}
\end{enumerate}
\end{Remark}

Let $X$ be a complex Banach space and $A$ be a bounded operator on $X$. For  $\lambda\in\mathbb{C}\backslash\sigma(A)\cup \sigma_{d}(T),$ where $\sigma_{d}(A)$ is the set of the Riesz points of $A$.
 We consider the operator
 \begin{align*}
 R_{B}(\lambda,A)=(A-\lambda \mathbb{I}_{X})(\mathbb{I}_{X}-P_{\{\lambda\}})+P_{\{\lambda\}}.
 \end{align*}
Thus, $R_{B}(\lambda,A)$ is invertible and
\begin{align*}
R_{B}^{-1}(\lambda,A)=(A-\lambda \mathbb{I}_{X}\mid_{N(P_{\{\lambda\}})})^{-1}(I-P_{\{\lambda\}})+P_{\{\lambda\}},
\end{align*}
where $P_{\{\lambda\}}$ is the Riesz projection complex associated to $\lambda.$ The Browder resolvent operator satisfies the Browder resolvent equation:
\begin{align*}
R_{B}^{-1}(\lambda,A)-R_{B}^{-1}(\mu,A)
&=(\lambda-\mu)R_{B}^{-1}(\lambda,A)R_{B}^{-1}(\mu,A) \\
&+R^{-1}(\lambda,A)[(A-(\lambda+1)\mathbb{I}_{X})P_{\{\lambda\}}-(A-(\mu+1)\mathbb{I}_{X})P_{\{\mu\}}]R^{-1}_{B}(\mu,A),
\end{align*}
for $\lambda$, $\mu$  $\in \rho(T) \cup \sigma_{d}(T)$.
\vskip 0.2 cm

Now, we give the Browder $S$-resolvent equation in quaternionic setting.
\begin{theorem}\label{th4.4}
Let $T\in\mathcal{B}(V^{R}_{\mathbb{H}})$ and let $s$, $p$ $\in\sigma_{d}^S(T)\cup\rho_S(T)$ with $p\not\in[s]$, then the left and right  Browder $S-$resolvent operators satisfies the following equation
\begin{align*}
 S^{-1}_{R,B}(s,T)S^{-1}_{L,B}(p,T)\mathcal{Q}_{s}(p)
                    &=[S^{-1}_{R,B}(s,T)- S^{-1}_{L,B}(p,T)] p+ \overline{s} [S^{-1}_{L,B}(p,T)-S^{-1}_{R,B}(s,T)]\\
                    &+[S^{-1}_{R,B}(s,T)(T-(p+1)\mathbb{I_{V^{R}_{\mathbb{H}}}})P_{[p]}-(T-(s+1)\mathbb{I_{V^{R}_{\mathbb{H}}}})P_{[s]} S^{-1}_{L,B}(p,T)]p\\
                    &+\overline{s}[(T-(s+1)\mathbb{I_{V^{R}_{\mathbb{H}}}})P_{[s]}S^{-1}_{L,B}(p,T)-S^{-1}_{R,B}(s,T)(T-(p+1)\mathbb{I_{V^{R}_{\mathbb{H}}}})P_{[p]}].
\end{align*}
\end{theorem}
\proof
Set
\begin{align*}
& \sigma_{B}(s,p,T):= S^{-1}_{R,B}(s,T)S^{-1}_{L,B}(p,T)(p^{2}-2Re(s)p+|s|^{2})\
\end{align*}
The left Browder $S$-resolvent equation implies
\begin{align*} S^{-1}_{L,B}(p,T)p=\mathbb{I_{V^{R}_{\mathbb{H}}}}+TS^{-1}_{L,B}(p,T)+(T-(p+1)\mathbb{I_{V^{R}_{\mathbb{H}}}})P_{[p]}. \end{align*}
In this way, we have
\begin{align*}
\sigma_{B}(s,p,T)&=S^{-1}_{R,B}(s,T)[\mathbb{I_{V^{R}_{\mathbb{H}}}}+ T S^{-1}_{L,B}(p,T)+(T-(p+1)\mathbb{I_{V^{R}_{\mathbb{H}}}})P_{[p]}]p\\
&-2Re(s)S^{-1}_{R,B}(s,T)[\mathbb{I_{V^{R}_{\mathbb{H}}}}+TS^{-1}_{L,B}(p,T)+(T-(p+1)\mathbb{I_{V^{R}_{\mathbb{H}}}})P_{[p]}]\\
&+|s|^{2}S^{-1}_{R,B}(s,T)S^{-1}_{L,B}(p,T)\\
&=S^{-1}_{R,B}(s,T)p+S^{-1}_{R,B}(s,T)T[\mathbb{I_{V^{R}_{\mathbb{H}}}}+TS^{-1}_{L,B}(p,T)+(T-(p+1)\mathbb{I_{V^{R}_{\mathbb{H}}}})P_{[p]}]\\
&+S^{-1}_{R,B}(s,T)(T-(p+1)\mathbb{I_{V^{R}_{\mathbb{H}}}})P_{[p]}p-2Re(s)S^{-1}_{R,B}(s,T)\\
&-2Re(s)S^{-1}_{R,B}(s,T)TS^{-1}_{L,B}(p,T)\\
&-2Re(s)S^{-1}_{R,B}(s,T)[(T-(p+1)\mathbb{I_{V^{R}_{\mathbb{H}}}})P_{[p]}]+|s|^{2}S^{-1}_{R,B}(s,T)S^{-1}_{L,B}(p,T).
 \end{align*}
The right Browder $S$-resolvent equation implies
\begin{align*}
S^{-1}_{R,B}(s,T)T=sS^{-1}_{R,B}(s,T)-(T-(s+1)\mathbb{I_{V^{R}_{\mathbb{H}}}})P_{[s]}-\mathbb{I_{V^{R}_{\mathbb{H}}}}.
\end{align*}
So, we get
\begin{align*}
\sigma_{B}(s,p,T)&=S^{-1}_{R,B}(s,T)p+ sS^{-1}_{R,B}(s,T)-(T-(s+1)\mathbb{I_{V^{R}_{\mathbb{H}}}})P_{[s]}-\mathbb{I_{V^{R}_{\mathbb{H}}}} \\
&+[ sS^{-1}_{R,B}(s,T)-(T-(s+1)\mathbb{I_{V^{R}_{\mathbb{H}}}})P_{[s]}-\mathbb{I_{V^{R}_{\mathbb{H}}}}] TS^{-1}_{L,B}(s,T)\\
&+[ sS^{-1}_{R,B}(s,T)-(T-(s+1)\mathbb{I_{V^{R}_{\mathbb{H}}}})P_{[s]}-\mathbb{I_{V^{R}_{\mathbb{H}}}}](T-(p+1)\mathbb{I_{V^{R}_{\mathbb{H}}}})P_{[p]} \\&+S^{-1}_{R,B}(s,T)(T-(p+1)\mathbb{I_{V^{R}_{\mathbb{H}}}})P_{[p]}p\\
&-2Re(s)S^{-1}_{R,B}(s,T)- 2Re(s)[ sS^{-1}_{R,B}(s,T)-(T-(s+1)\mathbb{I_{V^{R}_{\mathbb{H}}}})P_{[s]}-\mathbb{I_{V^{R}_{\mathbb{H}}}} ] S^{-1}_{L,B}(p,T)\\
&-2Re(s)S^{-1}_{R,B}(s,T)[T-(p+1)\mathbb{I_{V^{R}_{\mathbb{H}}}})P_{[p]}]+|s|^{2}S^{-1}_{R,B}(s,T)S^{-1}_{L,B}(p,T).
\end{align*}
If $s\not\in  [p],$ then $P_{[s]}P_{[p]}=0.$ In particular, we have
\begin{align*}(T-(s+1)\mathbb{I_{V^{R}_{\mathbb{H}}}})P_{[s]}(T-(p+1)\mathbb{I_{V^{R}_{\mathbb{H}}}})P_{[p]}= 0.\end{align*}
The left Browder $S$-resolvent equation implies
\begin{align*}
-\mathbb{I_{V^{R}_{\mathbb{H}}}} -T S^{-1}_{L,B}(p,T)=(T-(p+1)\mathbb{I_{V^{R}_{\mathbb{H}}}})P_{_{[s]}}-S^{-1}_{L,B}(p,T)p.
\end{align*}
Then, we obtain
\begin{align*}
\sigma_{B}(s,p,T)
&=S^{-1}_{R,B}(s,T)p+s S^{-1}_{R,B}(s,T)-(T-(s+1)\mathbb{I_{V^{R}_{\mathbb{H}}}})P_{[s]}-S^{-1}_{L,B}(p,T)p \\
&+ s[s S^{-1}_{R,B}(s,T)-(T-(s+1)\mathbb{I_{V^{R}_{\mathbb{H}}}})P_{[s]}-\mathbb{I_{V^{R}_{\mathbb{H}}}}]S^{-1}_{L,B}(p,T) \\
&-(T-(s+1)\mathbb{I_{V^{R}_{\mathbb{H}}}})P_{[s]}TS^{-1}_{L,B}(p,T) + s S^{-1}_{R,B}(s,T)(T-(p+1))P_{[p]}\\
&+S^{-1}_{R,B}(s,T)(T-(p+1))P_{[p]}p \\
&-2Re(s)S^{-1}_{R,B}(s,T)-2Re(s)sS^{-1}_{R,B}(s,T)S^{-1}_{L,B}(p,T)\\
&-2Re(s)[-(T-(s+1)\mathbb{I_{V^{R}_{\mathbb{H}}}})P_{[s]}S^{-1}_{L,B}(p,T)+S^{-1}_{R,B}(s,T)(T-(p+1)\mathbb{I_{V^{R}_{\mathbb{H}}}})P_{p}]\\
&+2Re(s)S^{-1}_{L,B}(p,T)+|s|^{2}S^{-1}_{R,B}(s,T)S^{-1}_{L,B}(p,T).\\
&=S^{-1}_{R,B}(s,T)p+sS^{-1}_{R,B}(s,T)-(T-(s+1)\mathbb{I_{V^{R}_{\mathbb{H}}}})P_{[s]}- S^{-1}_{L,B}(p,T)p\\
&+s^{2}S^{-1}_{R,B}(s,T)S_{L,B}^{-1}(p,T)-s(T-(s+1)\mathbb{I_{V^{R}_{\mathbb{H}}}})P_{[s]}S^{-1}_{L,B}(p,T)-sS^{-1}_{L,B}(p,T)\\
&-(T-(s+1)\mathbb{I_{V^{R}_{\mathbb{H}}}})P_{[s]}S^{-1}_{L,B}(p,T)p+(T-(s+1)\mathbb{I_{V^{R}_{\mathbb{H}}}})P_{[s]}\\
&+ sS^{-1}_{R,B}(s,T) (T-(p+1)\mathbb{I_{V^{R}_{\mathbb{H}}}})P_{[p]}+ S^{-1}_{R,B}(s,T)(T-(p+1)\mathbb{I_{V^{R}_{\mathbb{H}}}})P_{[p]}p\\
&-2Re(s)S^{-1}_{R,B}(s,T)-2Re(s)sS^{-1}_{R,B}(s,T)S^{-1}_{L,B}(p,T)\\
&-2Re(s)[-(T-(s+1)\mathbb{I_{V^{R}_{\mathbb{H}}}})P_{[s]}S^{-1}_{L,B}(p,T)+S^{-1}_{R,B}(s,T)(T-(p+1)\mathbb{I_{V^{R}_{\mathbb{H}}}})P_{[p]}]\\
&+2Re(s)S^{-1}_{L,B}(p,T)+|s|^{2}S^{-1}_{R,B}(s,T)S^{-1}_{L,B}(p,T).\\
&=[S^{-1}_{R,B}(s,T)-S^{-1}_{L,B}(p,T)]p+(s^{2}-2Re(s)s+|s|^{2})S^{-1}_{R,B}(s,T)S^{-1}_{L,B}(p,T)\\
&+[S^{-1}_{R,B}(s,T)(T-(p+1)\mathbb{I_{V^{R}_{\mathbb{H}}}})P_{[p]}-(T-(s+1)\mathbb{I_{V^{R}_{\mathbb{H}}}})P_{[s]}S^{-1}_{L,B}(p,T)]p-\overline{s}S^{-1}_{R,B}(s,T)\\
&+  \overline{s}[(T-(s+1)\mathbb{I_{V^{R}_{\mathbb{H}}}})P_{[s]}S^{-1}_{L,B}(p,T)-S^{-1}_{R,B}(s,T)(T-(p+1)\mathbb{I_{V^{R}_{\mathbb{H}}}})P_{[p]}]+\overline{s}S^{-1}_{L,B}(p,T).
\end{align*}
On the other hand
\begin{align*}s^2+2{\rm Re}(s)s+\vert s\vert^2=0.\end{align*}
We conclude that
\begin{align*}
 \sigma_{B}(s,p,T)
                    &=[S^{-1}_{R,B}(s,T)- S^{-1}_{L,B}(p,T)] p+ \overline{s} [S^{-1}_{L,B}(p,T)-S^{-1}_{R,B}(s,T)]\\
                    &+[S^{-1}_{R,B}(s,T)(T-(p+1)\mathbb{I_{V^{R}_{\mathbb{H}}}})P_{[p]}-(T-(s+1)\mathbb{I_{V^{R}_{\mathbb{H}}}})P_{[s]} S^{-1}_{L,B}(p,T)]p\\
                    &+\overline{s}[(T-(s+1)\mathbb{I_{V^{R}_{\mathbb{H}}}})P_{[s]}S^{-1}_{L,B}(p,T)-S^{-1}_{R,B}(s,T)(T-(p+1)\mathbb{I_{V^{R}_{\mathbb{H}}}})P_{[p]}].
\end{align*}\qed
\begin{Remark}\label{rem4.5}
\begin{enumerate}
\item If $s$, $p$ $\in\rho_{S}(T),$ then $P_{[s]}=P_{[p]}=0.$
Hence, we find the $S$-resolvent equation, see \cite[theorem 3.1.15]{FJDP}:\\
$S^{-1}_{R}(s,T)S^{-1}_{L}(p,T)$
\begin{align*}
=[(S^{-1}_{R}(s,T)-S^{-1}_{L}(p,T))p-\overline{s}[S^{-1}_{R}(s,T)-S^{-1}_{L}(p,T)](p^{2}-2Re(s)p+|s|^{2})^{-1}.
\end{align*}
\item Let's test the Browder $S$-resolvent equation in the commutative case,
  if Tq=qT for all $q\in\mathbb{H},$ then for $s$, $q$ $\in \rho_{S}(T),$ we have
\begin{align*}
M(s,p)&:=[S^{-1}_{R,B}(s,T)(T-(p+1)\mathbb{I_{V^{R}_{\mathbb{H}}}})P_{[p]}-(T-(s+1)\mathbb{I_{V^{R}_{\mathbb{H}}}})P_{[s]}S_{R,B}^{-1}(p,T)]p\\
&+\overline{s}[(T-(s+1)\mathbb{I_{V^{R}_{\mathbb{H}}}})P_{[s]} S^{-1}_{L,B}(p,T)- S^{-1}_{R,B}(s,T)(T-(p+1)\mathbb{I_{V^{R}_{\mathbb{H}}}})P_{[p]}]\\
&=p[S^{-1}_{R,B}(s,T)(T-(p+1)\mathbb{I_{V^{R}_{\mathbb{H}}}})P_{[p]}-S^{-1}_{L,B}(p,T)(T-(s+1)\mathbb{I_{V^{R}_{\mathbb{H}}}})P_{[s]}]\\
&+\overline{s}[S^{-1}_{L,B}(p,T)(T-(s+1)\mathbb{I_{V^{R}_{\mathbb{H}}}})P_{[s]}-S^{-1}_{R,B}(s,T)(T-(p+1)\mathbb{I_{V^{R}_{\mathbb{H}}}})P_{[p]}]\\
&=(p-\overline{s})[S^{-1}_{R,B}(s,T)(T-(p+1)\mathbb{I_{V^{R}_{\mathbb{H}}}})P_{[p]}-S^{-1}_{L,B}(p,T)(T-(s+1)\mathbb{I_{V^{R}_{\mathbb{H}}}})P_{[s]}].
\end{align*}
\end{enumerate}
\end{Remark}
\noindent In particular in the complex case, if $T\in\mathcal{B}(V_{\mathbb{C}})$ then,
\begin{align*}S_{R,B}^{-1}(p,T)=R_{B}^{-1}(p,T)\mbox{ and }S_{L,B}^{-1}(s,T)= R_{B}^{-1}(s,T). \end{align*}
Therefore, we obtain\\
$M(s,p)$
\begin{align*}
&=(p-\overline{s})R_{B}^{-1}(s,T)[(T-(p+1)\mathbb{I_{V^{R}_{\mathbb{H}}}})P_{[p]}- R_{B}(s,T)(T-(s+1)\mathbb{I_{V^{R}_{\mathbb{H}}}})P_{[s]}R^{-1}_{B}(p,T)]\\
&=(p-\overline{s})R^{-1}_{B}(s,T)[(T-(p+1)\mathbb{I_{V^{R}_{\mathbb{H}}}}P_{[s]}R_{B}(p,T)-R_{B}(s,T)(T-(s+1)\mathbb{I_{V^{R}_{\mathbb{H}}}})P_{[s]}]R^{-1}_{B}(p,T).
\end{align*}
Hence,
$$P_{[p]}R_{B}(p,T)=P_{[p]}\ \mbox{and}\ R_{B(s,T)}P_{[s]}=P_{[s]}.$$
So, we get
\begin{align*}
M(s,p):=(p-\overline{s})R^{-1}_{B}[(s,T)[(T-(p+1)\mathbb{I_{V^{R}_{\mathbb{H}}}})P_{[p]}-(T-(s+1)\mathbb{I_{V^{R}_{\mathbb{H}}}})P_{[s]}]R^{-1}_{B}(p,T).
\end{align*}
Thus, we obtain the classic Browder resolvent equation in complex case.

\end{document}